\documentclass[letterpaper, preprint, paper,11pt]{AAS}

\usepackage{mathtools}
\usepackage{gensymb}
\usepackage{float}
\usepackage{graphicx}
\usepackage{xcolor}
\usepackage{tabularx}
\usepackage{adjustbox}

\usepackage{bm}
\usepackage{amsmath}
\usepackage{subfigure}
\usepackage{comment}
\usepackage{overcite}
\usepackage{footnpag}			      	
\usepackage{amssymb}
\usepackage{subfigure}   
\usepackage{multicol}
\usepackage{multirow}
\usepackage{booktabs}
\usepackage{pstool}             
\usepackage[colorlinks]{hyperref}

\usepackage{amsthm}

\PaperNumber{25-123}

\makeatletter
\newcommand{\ostar}{\mathbin{\mathpalette\make@circled*}}
\newcommand{\make@circled}[2]{%
  \ooalign{$\m@th#1\smallbigcirc{#1}$\cr\hidewidth$\m@th#1#2$\hidewidth\cr}%
}
\newcommand{\smallbigcirc}[1]{%
  \vcenter{\hbox{\scalebox{0.77778}{$\m@th#1\bigcirc$}}}%
}
\makeatother
\usepackage{subfigure}

\begin{document}

\title{Advances in Cislunar Periodic Solutions via Taylor Polynomial Maps}

\author{ Mohammed Atallah\thanks{PhD Student, Department of Aerospace Engineering, Iowa State University, IA 50011, USA. email: matallah@iastate.edu} 
\ and Simone Servadio\thanks{Assistant Professor, Department of Aerospace Engineering, Iowa State University, IA 50011, USA. email: servadio@iastate.edu}
}

 \maketitle
\begin{abstract}
{In this paper, novel approaches are developed to explore the dynamics of motion in periodic orbits near libration points in cislunar space using the Differential Algebra (DA) framework. The Circular Restricted Three-Body Problem (CR3BP) models the motion, with initial states derived numerically via differential correction. Periodic orbit families are computed using the Pseudo-Arclength Continuation (PAC) method and fitted. Two newly developed polynomial regression models (PRMs) express initial states as functions of predefined parameters and are used in the DA framework to evaluate propagated states. The initial states, expressed via PRM, are propagated in the DA framework using the fourth-order Runge-Kutta (RK4) method. The resultant polynomials of both PRM and DA are employed to develop a control law that shows significantly reduced control effort compared to the traditional tracking control law, demonstrating their potential for cislunar space applications, particularly those requiring computationally inexpensive low-energy transfers.}

\end{abstract}
 
\section{Introduction}
With the increasing potential for deep-space exploration missions, there is significant interest in cislunar space due to its role in designing low-energy trajectories [\citen{koon2000dynamical}]. However, this region is recognized as a chaotic system because of its multi-body gravitational environment. These characteristics have drawn attention to the bounded motion represented by the periodic and quasi-periodic orbits near libration points~[\citen{wilmer2024preliminary}]. In the past few decades, numerous studies have investigated motion in cislunar space. The Circular Restricted Three-Body Problem (CR3BP) is one of the simplified mathematical models commonly used to find solutions for bounded trajectories. This model is linearized around a libration point to obtain the trajectory of the periodic orbit, which introduces inaccuracy. This inaccuracy is compensated using a high-order differential correction scheme that obtains a more accurate trajectory [\citen{richardson1980halo}]. However, even the differential correction scheme cannot provide a high-fidelity solution for the trajectory due to unmodeled dynamics and external disturbances. Therefore, there is a need for a real-time correction scheme that retains the satellite in a periodic orbit by leveraging sensor measurements.

In recent years, several missions have utilized periodic orbits near the libration points in cislunar space. For instance, the first stationkeeping operations around \(\mathcal{L}_1\) and \(\mathcal{L}_2\) in cislunar space were performed by the ARTEMIS mission [\citen{woodard2009artemis}]. In light of these advancements, the Lunar Orbiter Platform-Gateway (LOP-G) is one of the largest international cooperative space programs, aiming to assemble a space station around the Moon [\citen{merri2018lunar}]. Additionally, NASA's journey to Mars will utilize cislunar space to conduct advanced operations [\citen{national2016nasa}]. Moreover, the next decade will witness over thirty missions being launched in the cislunar region [\citen{baker2024comprehensive}]. These missions require staging locations to conduct various activities, where the periodic orbits near the libration points are being investigated as potential choices for that role [\citen{whitley2016options}].

The CR3BP is the most commonly used framework for transfers in cislunar space, where approximate trajectory solutions can be obtained analytically [\citen{richardson1980analytic}]. In [\citen{pritchett2018impulsive}], a methodology for conducting low-energy transfers between periodic orbits is developed using CR3BP. In [\citen{singh2020low}], the Halo orbit (HO) family near \(\mathcal{L}_1\) is utilized to design a low-thrust transfer of a small spacecraft to a low-altitude lunar orbit. In [\citen{davis2017stationkeeping}], stationkeeping in cislunar space is investigated, with the CR3BP being employed to generate the periodic orbit families, while the n-body dynamical model simulates higher-fidelity trajectories. In [\citen{van2016tadpole}], the Pseudo-Arclength Continuation (PAC) method is developed to compute the members across each periodic orbit family based on the CR3BP. In [\citen{wilmer2021cislunar}], the benefits of \(\mathcal{L}_1\) and \(\mathcal{L}_2\) HOs for orbit maintenance are investigated, where the CR3BP is used to model the constellations. In [\citen{fay2024investigation}], search and rescue operations are investigated, and the response times are compared for rescuer spacecraft located in distant retrograde orbits and \(\mathcal{L}_1\)/\(\mathcal{L}_2\) Lyapunov orbit (LO) families. More accurate trajectory solutions can be obtained using the bicircular restricted four-body problem (BCR4BP), as presented in [\citen{negri2020generalizing, oshima2022multiple, wilmer2021lagrangian}]; however, these solutions are more computationally expensive and cannot be obtained analytically. Therefore, this study employs the CR3BP to develop a methodology for representing motion in a periodic orbit family leveraging Differential Algebra (DA).

DA is a computationally efficient tool based on Taylor expansion, that can be employed to represent differentiable and continuous dynamic models as high-order polynomials~[\citen{chao2023handbook, hawkes1999modern}]. Several tools are supplied in the DA framework to obtain the derivatives and integrals of the models in low-level computation environments, such as FORTRAN~[\citen{berz1987differential}], and C/C++~[\citen{massari2018differential}]. In addition, DA has been proven to be a reliable tool for numerical integration of Ordinary Differential Equations (ODE) carried out by an arbitrary integration scheme. Several applications have leveraged DA framework, such as describing beam dynamics~[\citen{berz1988differential}], and high-order nonlinear filtering~[\citen{valli2014nonlinear, cavenago2018based, servadio2021differential}]. The DA tool is fundamentally based on expressing a continuous differentiable function as an infinite series expanded at a predetermined operating point [\citen{servadio2022maximum}]. For a small deviation of this operating point, the series returns a precise value of the function using a finite number of the terms. This introduces the concept of Truncated Power Series (TPS), which is computationally reliable and can be used for applying arithmetic and calculus operations. In this study, the DA is employed to represent the initial states of periodic orbit families as functions of predetermined parameters, then these TPS are propagated according to ODE of the CR3BP using fourth-order Runge-Kutta scheme.

This paper aims to investigate and analyze the periodic orbit families near libration points in cislunar space within the framework of DA. First, the general translational motion in cislunar space is expressed using the CR3BP. Then, the initial states of an arbitrary periodic orbit in a given family are obtained using the linearized model around the nearest libration point. These approximate states are refined using a high-order differential correction scheme to obtain more accurate ones. Next, the members across the family are computed using the Pseudo-Arclength Continuation (PAC) method. After that, these members are employed to fit a Polynomial Regression Model (PRM) using the Least-Squares Error (LSE) method. The resultant polynomials of the periodic orbit initial states are then propagated to specific times using the fourth-order Runge-Kutta scheme within the DA framework. The first propagation process uses absolute time, while the second process uses normalized time, in which each orbit in the family is propagated for a fractional amount of its period. Finally, numerical simulations are conducted to demonstrate the reliability and accuracy of DA in representing different periodic orbit families near libration points in cislunar space. Additionally, a Proportional-Derivative (PD) control law is developed using the proposed method and compared to the traditional tracking control law to demonstrate the optimality of the proposed approach.

The rest of the paper is organized as follows: Section~2 presents the mathematical model of the CR3BP. Section~3 introduces the basics of DA. Section~4 shows the evaluation of the periodic orbits near $\mathcal{L}_1$ and $\mathcal{L}_2$ and the application of DA in that process. Section~5 presents and discusses the results of the numerical simulations and demonstrates the applicability of the proposed methodology. Section~6 concludes the paper.

\section{Mathematical Model of the CR3BP}
The translational motion in cislunar space can be approximated by an autonomous dynamic model by applying the following assumptions: 
\begin{enumerate}
    \item The Earth and the Moon are treated as mass points.
    \item The Moon moves in a circular orbit around the Earth.
    \item The gravity of the Earth and the Moon is the only source of force influencing the motion, while all other perturbations are neglected.
\end{enumerate}

Conventionally, the parameters and states in the CR3BP model are dimensionless, and the motion is {expressed in rotating coordinates} centered at the barycenter of the Earth-Moon system. The \(X\) axis is in the direction of the vector between the Earth and the Moon. Equation~\eqref{eq:CR3BP} presents the mathematical model of the CR3BP according to the aforementioned assumptions and conventions.
\begin{equation}
\begin{aligned}
& \ddot{x}=2 \dot{y}+x-\dfrac{(1-\mu)(x+\mu)}{r_1^3}-\dfrac{\mu[x-(1-\mu)]}{r_2^3} \\
& \ddot{y}=-2 \dot{x}+y-\dfrac{(1-\mu) y}{r_1^3}-\dfrac{\mu y}{r_2^3} \\
& \ddot{z}=-\dfrac{(1-\mu) z}{r_1^3}-\dfrac{\mu z}{r_2^3}
\end{aligned}
\label{eq:CR3BP}  
\end{equation}
Here, \(x\), \(y\), and \(z\) denote the components of the dimensionless position vector of the satellite, where \(x\) points to the Moon, \(y\) is in the direction of the relative motion of the Moon with respect to the Earth, and \(z\) completes the set according to the right-hand rule. {\(\mu = 0.01215\)} is the dimensionless mass of the Moon. \(r_1\) and \(r_2\) are the relative distances between the satellite and the Earth, and the satellite and the Moon, respectively.

\section{Periodic Orbits Near $\mathcal{L}_1$ and $\mathcal{L}_2$}
In the CR3BP model, there are five equilibrium points, known as libration points, where the gravitational forces exerted by the Earth and the Moon on a satellite are balanced. The first two points, \(\mathcal{L}_1\) and \(\mathcal{L}_2\), have special characteristics due to their symmetry relative to the Moon. There are two common families of periodic orbits near these points: LOs, which exist in two-dimensional space~[\citen{henon1969numerical}], and HOs, which exist in three-dimensional space~[\citen{breakwell1979halo}]. Computing these periodic orbits requires a series of iterative steps due to the chaotic behavior of the CR3BP and the absence of an analytical solution for the model, as follows:

\begin{enumerate}
    \item The mathematical equations are linearized at the libration point, and approximate initial states are computed using this linear model.
    \item An iterative high-order differential correction scheme is employed to determine the period and refine the initial states of the orbit.
    \item The computed member of the orbit family is used to generate other members in the family using the PAC method.
\end{enumerate}

\subsection{Linearized Equations of Motion}
The detailed steps of linearizing the model are found in~[\citen{richardson1980analytic, richardson1980halo}]. The resultant linear model is derived as follows:
\begin{equation}
\begin{aligned}
&\ddot{x} = 2 \dot{y}+\left(1+2 c_2\right) x\\
&\ddot{y} = -2 \dot{y}-\left(c_2-1\right) y \\
&\ddot{z} = -c_2 z
\end{aligned}
\label{eq:linModel}
\end{equation}
where  
\begin{equation}
c_n=\dfrac{1}{\gamma_L^3}\left[( \pm 1)^n \mu+(-1)^n \dfrac{(1-\mu) \gamma_L^{n+1}}{\left(1 \mp \gamma_L\right)^{n+1}}\right], \quad\left(L_1 \text { or } L_2\right)
\end{equation}
Here, $\mu_E=G M_E$, $G$ is the gravitational constant, $M_E$ is the Earth mass, {$\gamma_L=r_E / a$, $r_E$} is the Earth mean radius, and $a$ is the astronomical unit.

\subsection{Differential Correction for Computing Initial States}
The linearized model in Equation~\eqref{eq:linModel} has a closed-form analytical solution, from which an analytical formula for the initial states of periodic orbits can be derived. However, the accuracy of these approximate initial states is insufficient due to the chaotic behavior of the system, {and the high non-linearity of the real-time system.} Therefore, further correction is required to achieve the desired accuracy. The differential correction method, commonly used for this purpose, implements an iterative algorithm using the nonlinear model to modify the initial states. The differential correction scheme used in this study was first proposed in~[\citen{richardson1980halo}]. The procedure of this scheme is as follows: First, the initial states of the periodic orbit are obtained using the closed-form analytical solution of Equation~\eqref{eq:linModel}. These states are then propagated using the CR3BP model until they intersect the x-z plane for HOs or the y-axis for LOs. {Due to the symmetry of the periodic orbit, the states \(\dot{x}\) and \(\dot{z}\) must equal zero at the intersection point.} Next, the state transition matrix is evaluated at this half-period. Using the states and the state transition matrix at the half-period point, the next iteration of the corrected initial states is computed as follows:
\begin{equation}
\left[\begin{array}{c}
\Delta x_0 \\
\Delta \dot{y}_0 \\
\Delta T_{1 / 2}
\end{array}\right]=-\Phi^{-1}\left[\begin{array}{c}
\dot{x}\left(x_0, \dot{y}_0, T_{1 / 2}\right) \\
\dot{z}\left(x_0, \dot{y}_0, T_{1 / 2}\right) \\
y\left(x_0, \dot{y}_0, T_{1 / 2}\right)
\end{array}\right]
\end{equation}
where $\Phi$ is the matrix of partials that is defined as follows:
\begin{equation}
\Phi=\left\{\begin{array}{ccc}
\frac{\partial \dot{x}}{\partial x_0} & \frac{\partial \dot{x}}{\partial \dot{y}_0} & \frac{\partial \dot{x}}{\partial T_{1 / 2}} \\
\frac{\partial \dot{z}}{\partial x_0} & \frac{\partial \dot{z}}{\partial \dot{y}_0} & \frac{\partial \dot{z}}{\partial T_{1 / 2}} \\
\frac{\partial y}{\partial x_0} & \frac{\partial y}{\partial \dot{y}_0} & \frac{\partial y}{\partial T_{1 / 2}}
\end{array}\right\}_{t=T_{1 / 2}}
\end{equation}
This process repeats until the desired state error is achieved.

\subsection{Computing Members across the Orbit Families}
The differential correction scheme is used obtain the initial states of a single periodic orbit. It might be used to obtain the other members of the periodic orbit family, however, it has a limited range of these members and it is not the most computationally efficient method for that purpose~[\citen{servadio2022dynamics}]. The {PAC} method is one of the most reliable methods that is developed to compute a wide range of members across the family using a single predetermined member~[\citen{van2016tadpole, patel2023low}]. It computes the members with a step \(\Delta s\), that is predetermined depending on the desired number of members, in the tangent direction to the solution manifold. {In this paper,} the method is employed for both Lyapunov and HOs, {though it can be used for all periodic and quasi-periodic orbit families.} 

Assume \(\mathbf{x}_{i}\) is the initial state vector of an arbitrary member that satisfies the constraints \(F\left(\mathbf{x}_{i}\right) = 0\) of the family. {In case of LO, the constraints are \(y|_{t=T/2} = \dot{x}|_{t=T/2} = 0\), while \(\dot{z}|_{t=T/2} = 0\) is added for HOs.} Shifting this member by a step size \(\Delta s\) yields the next member \(\mathbf{x}_{i+1}\), which also satisfied the constraints \(F\left(\mathbf{x}_{i+1}\right) = 0\). To guarantee that the step size equals \(\Delta s\) in the tangent direction, an additional constraint is added as follows: 
\begin{equation}
G\left(\mathbf{x}_{i+1}\right)=\left[\begin{array}{c}
F\left(\mathbf{x}_{i+1}\right) \\
\left(\mathbf{x}_{i+1}-\mathbf{x}_{i}\right)^T \Delta \mathbf{x}_{i}-\Delta s
\end{array}\right]=\mathbf{0}
\label{eq:G_cons}
\end{equation}
where {\(G(\cdot)\)} is the augmented constraints, and \(\Delta \mathbf{x}_{i}\) is the null vector of the Jacobian matrix for \(\mathbf{x}_{i}\), which is defined as $
\Delta \mathbf{x}_{i}=\mathfrak{N}\left(D_ F\left(\mathbf{x}_{i}\right)\right)
$. {Here, \(\mathfrak{N}\left(\cdot\right)\) denotes the null vector, and \(D_{(\cdot)}\) is the Jacobian matrix of \((\cdot)\).} The new Jacobian matrix of the augmented constraints is defined as follows:
\begin{equation}
D_G\left(\mathbf{x}_{i+1}\right)=\left[\begin{array}{c}
D_F\left(\mathbf{x}_{i+1}\right) \\
\Delta \mathbf{x}_{i}^{T}
\end{array}\right]
\end{equation}
In order to obtain \(\mathbf{x}_{i+1}\) that satisfies the constraints in Equation~\eqref{eq:G_cons}, an initial guess of the solution \({}^{0}\mathbf{x}_{i+1}\) is selected. Then, Newton's method is employed as follows:
\begin{equation}
^{k+1}\mathbf{x}_{i+1}={}^{k}\mathbf{x}_{i+1}-\left[D_G\left({}^{k}\mathbf{x}_{i+1}\right)\right]^{-1} G\left({}^{k}\mathbf{x}_{i+1}\right)
\end{equation}
where \({}^{k}\mathbf{x}_{i+1}\) is the \(k^\text{th}\) iteration. This process is repeated iteratively until the constraints in Equation~\eqref{eq:G_cons} are satisfied within a certain tolerance.

Starting from an arbitrary member in the family, the other members can be computed in both directions \(\pm \Delta \mathbf{x} \). In this study, the PAC method is employed to compute the members across LO families, HO, and NRHO families near \(\mathcal{L}_1\) and \(\mathcal{L}_2\).

\subsubsection{LO Families}
\begin{figure}[H]
    \centering
    \subfigure[Orbits near L1]{%
        \includegraphics[width=0.45\textwidth]{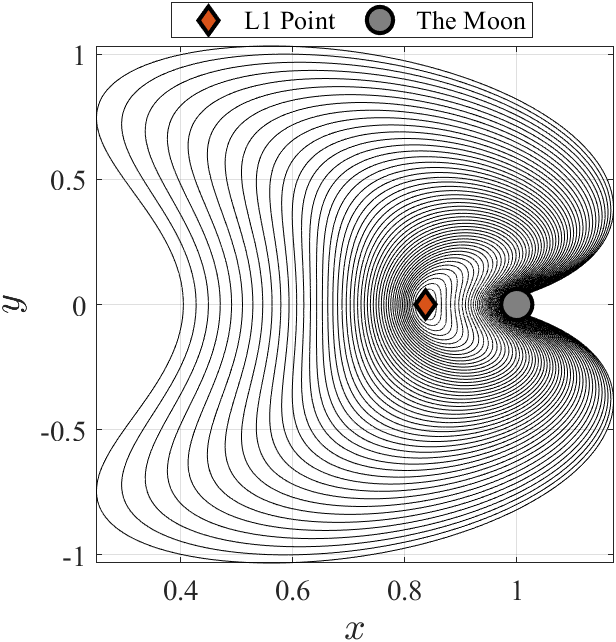}
    }
    \hfill
    \subfigure[Orbits near L2]{%
        \includegraphics[width=0.45\textwidth]{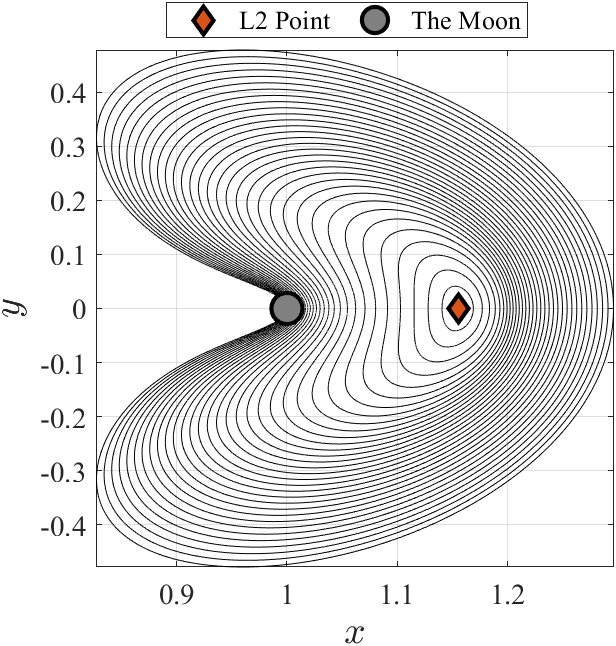}
    }
    \caption{The members of the LO Families.}
    \label{fig:Lyap}
\end{figure}

Figure~\ref{fig:Lyap} shows several LOs near \(\mathcal{L}_1\) and \(\mathcal{L}_2\). A key advantage of this method is its ability to obtain orbits that are close to the Moon.

\subsubsection{HO and NRHO Families} 
\begin{figure}[H]
    \centering
    \subfigure[Orbits near L1]{%
        \includegraphics[width=0.45\textwidth]{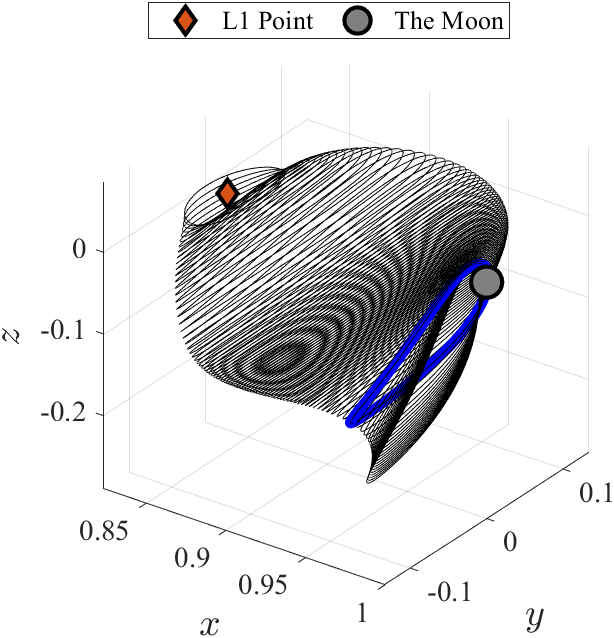}
        \label{fig:Halo_L1}
    }
    \hfill
    \subfigure[Orbits near L2]{%
        \includegraphics[width=0.45\textwidth]{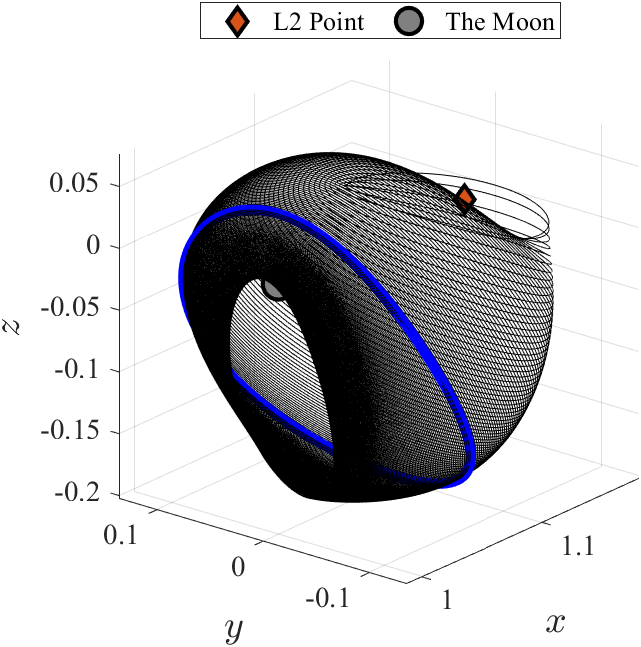}
        \label{fig:Halo_L2}
    }
    \caption{The members of the HO families.}
    \label{fig:Halo}
\end{figure}

\begin{figure}[H]
    \centering
    \includegraphics[width=0.6\linewidth]{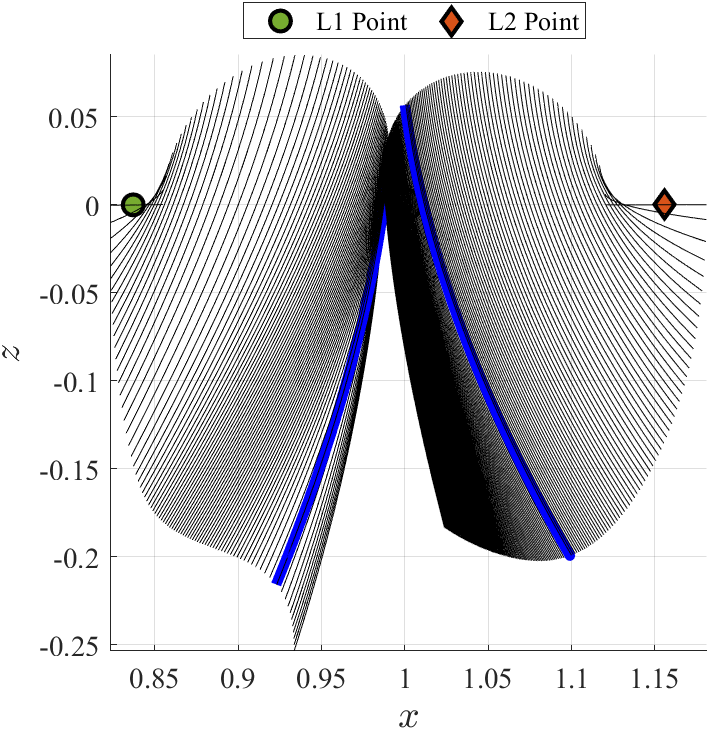}
    \caption{HO and NRHO families near L1 and L2.}
    \label{fig:aLLHalo_Orbits}
\end{figure}

Similarly, Figures~\ref{fig:Halo} and~\ref{fig:aLLHalo_Orbits} show numerous HOs and NRHOs near \(\mathcal{L}_1\) and \(\mathcal{L}_2\). Specifically, Figure~\ref{fig:Halo_L1} displays the \(\mathcal{L}_1\) families, separated by the bold blue orbit, while Figure~\ref{fig:Halo_L2} shows the \(\mathcal{L}_2\) families.

\section{Overview of Differential Algebra}
\begin{figure}[H]
    \centering
    \subfigure[Floating point mapping.]{%
        \includegraphics[width=0.42\textwidth]{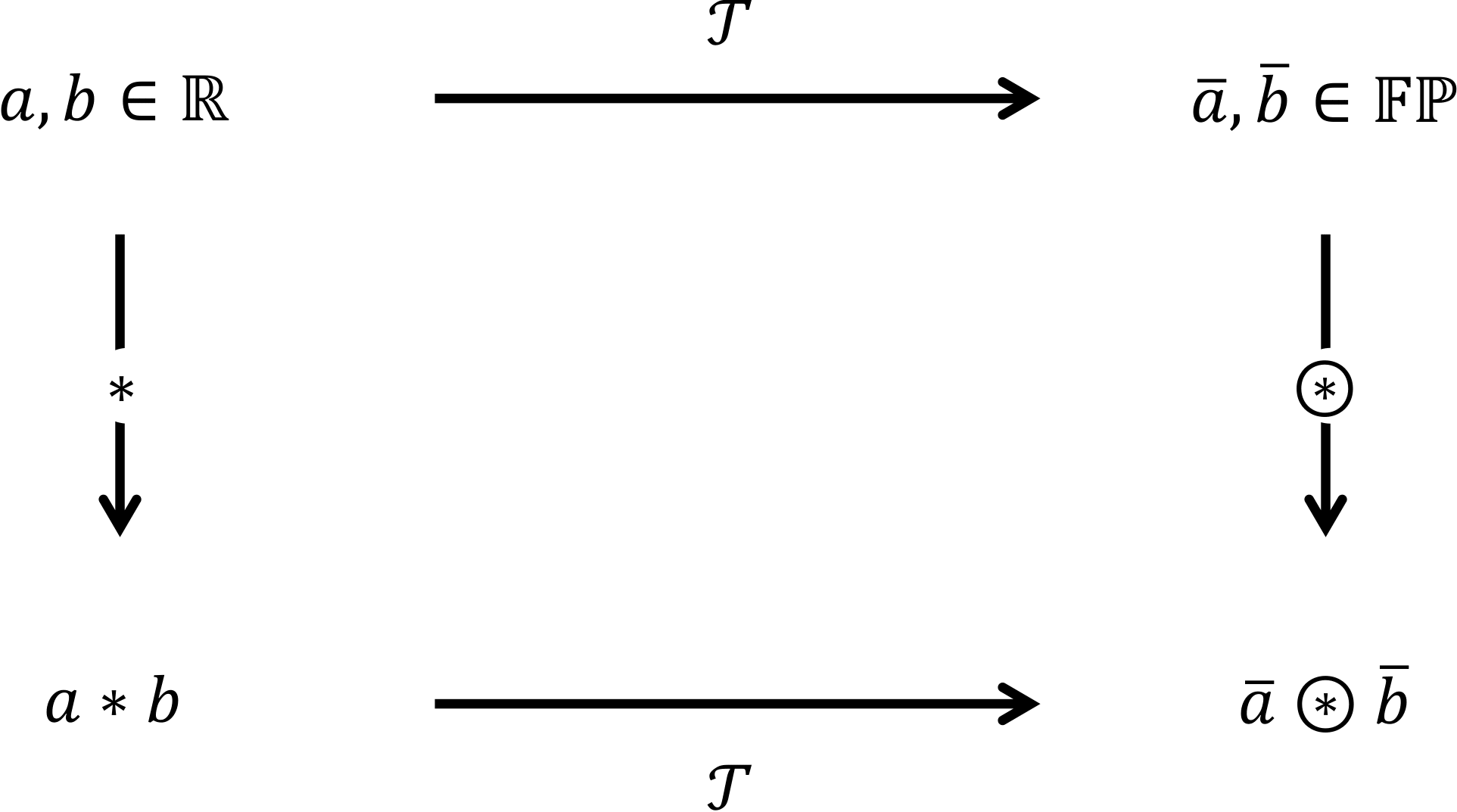}
        \label{fig:FPN}
    }
    \hfill
    \subfigure[DA mapping.]{%
        \includegraphics[width=0.5\textwidth]{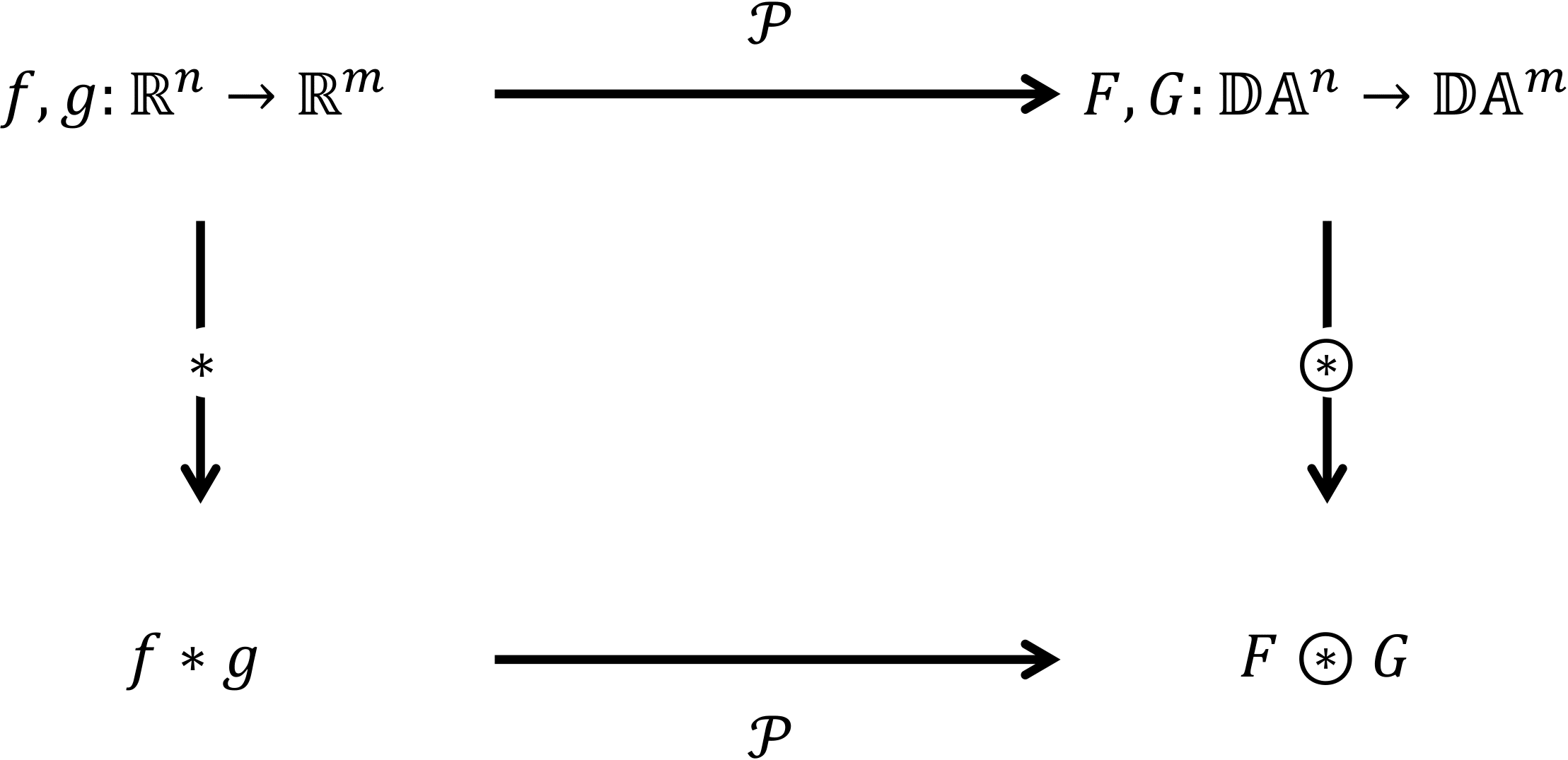}
        \label{fig:DAF}
    }
    \caption{Implementing operations on real numbers and functions by transforming them to floating point numbers or DA series.}
    \label{fig:Mapping}
\end{figure}

The basic idea of treating numbers and implementing various operations on them using computers is to represent these numbers with a finite amount of information. Generally, numbers can be irrational and ideally represented by infinite digits, which makes it impractical for computers to handle their ideal representation. Therefore, only a finite amount of relevant information is extracted. These approximations are known as floating-point numbers. This approximated form of the numbers allows operations on real numbers by transforming these numbers into floating-point numbers and implementing the operations as depicted in Figure~\ref{fig:FPN}. {Here, \(a\) and \(b\) are real numbers, \(\bar{a}\) and \(\bar{b}\) are floating-point numbers, and \(\ostar\) denotes an arbitrary operation.}

In a similar manner, the DA technique extends the concept of floating-point numbers to encompass differentiable functions [\citen{servadio2020recursive}]. According to the Taylor expansion, any differentiable function at a certain point can be represented by an infinite series expanded at that point. This brings an analogy between real numbers and differentiable functions, in which both are represented by an infinite amount of information (i.e., digits of real numbers and series coefficients of differentiable functions). Similar to floating-point numbers, DA extracts a finite number of terms to represent the function in an approximated way that can be handled by computers. This approximation is used to implement various operations on these functions in a computationally efficient manner. Figure~\ref{fig:DAF} depicts the equivalent DA approximation of functions in computer environments~[\citen{servadio2020nonlinear}]. {Here, \(f\) and \(g\) are differentiable functions, while \(F\) and \(G\) are finite series that represent these functions in the DA framework, defined by their coefficients.}

\begin{figure}[H]
    \centering
    \includegraphics[width=\linewidth]{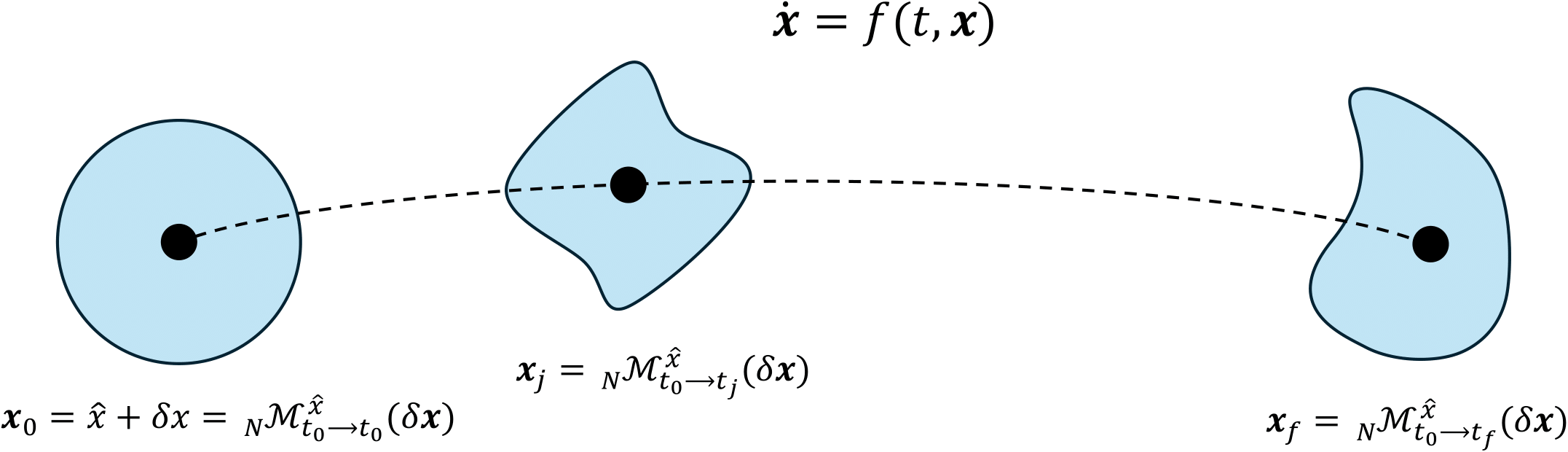}
    \caption{Propagating the states in the DA framework.}
    \label{fig:mappingDemo}
\end{figure}

{For any differentiable function \(\mathbf{y} = f\left(\mathbf{x}\right)\), the DA mapping is represented as follows:} 
\begin{equation}
    \mathbf{y}\left(\delta \mathbf{x}\right) = {}_{N}\mathcal{F}^{\hat{\mathbf{x}}}\left(\delta\mathbf{x}\right)
    \label{eq:mappingDemo}
\end{equation}
where \(\hat{\mathbf{x}}\) denotes the operating point of the series expansion, \(\delta \mathbf{x}\) denotes the deviation of the \(\mathbf{x}\) defined as \(\delta \mathbf{x} = \mathbf{x} - \hat{\mathbf{x}}\), and \(N\) is the highest order of the series with nonzero coefficient. This representation can be used to derive a highly accurate approximation of the solution for a dynamical system [\citen{servadio2022nonlinear}]. For a given dynamic system \(\dot{\mathbf{x}} = f\left(t, \mathbf{x}\right)\), the DA framework can be employed to evaluate the states at a certain time \(t_j\), as depicted in Figure~\ref{fig:mappingDemo}.
{Here, \(\mathcal{M}^{\hat{\mathbf{x}}}_{t_0 \rightarrow t_j}\left(\delta\mathbf{x}\right)\) denotes the State Transition Propagation Matrix (STPM) of the states propagated from the initial time \(t_0\) to a given time \(t_j\), expressed as a function of the deviation \(\delta \mathbf{x}\) with respect to the operating point of the expansion \(\hat{\mathbf{x}}\). This approach allows for the propagation of the neighborhood of a given state \(\mathbf{x}\) to multiple times by propagating the series instead of propagating each point individually [\citen{servadio2024likelihood}]. This method is computationally efficient in any application that requires computing multiple states at different times. Additionally, the accuracy can be balanced with computation time by tuning the order \(N\) of the series.}

\section{Polynomial Regression Model}
The computed members in each periodic orbit family are used to construct PRMs. Two different approaches are implemented: The first is the global PRM, where the domain of the predefined parameter \(\kappa\) is divided into multiple regions. The mean of each region is selected as the operating point \(\hat{\kappa}\), and the polynomial for that region is fitted. The deviations of the predefined parameter \(\delta \kappa\) for the members act as the query points, while the initial states of the members serve as the fitted values at these query points. The second approach is the local PRM, which uses the parameter \(\kappa\) at the designed periodic orbit as the operating point for polynomial fitting. However, this approach only uses the neighbors of the designed member to fit the polynomial, in order to avoid fitting issues.

\subsection{Global Polynomial Regression Model}
\begin{figure}[H]
    \centering
    \includegraphics[width=\linewidth]{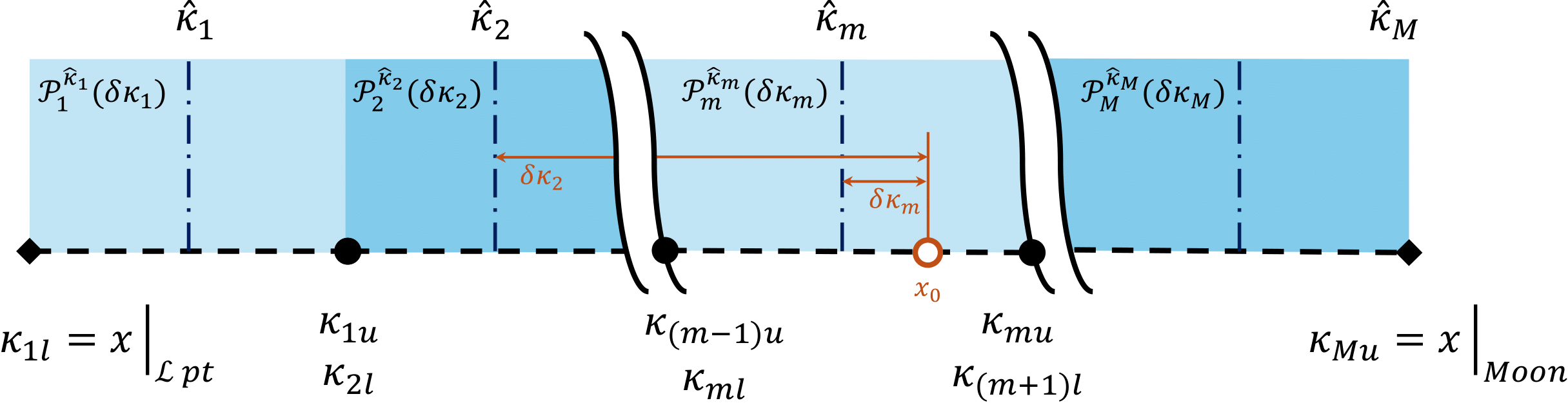}
    \caption{The regions of the global PRM that are constructed along the predefined parameter \(\kappa\), and defined by their upper and lower bounds, \(\kappa_u\) and \(\kappa_l\), as well as operating points, \(\hat{\kappa}\). }   \label{fig:regions_poly}
\end{figure}

In this study, the \(x\) component of the initial states of the members is used as the parameter \(\kappa\) for both Lyapunov and Halo families. The \(x\) component is chosen because it is unique for each member, unlike the \(y\) component, which is always zero, and the \(z\) component, which is zero in Lyapunov families and not unique in Halo families. Figure~\ref{fig:regions_poly} illustrates the concept of dividing the domain into multiple regions and fitting a polynomial in each, where \(\kappa_{il}\) and \(\kappa_{iu}\) represent the lower and upper bounds of the \(i\)th region {of the domain}, respectively, while \(\mathcal{P}^{\hat{\kappa}_i}_i\) is the polynomial of the members in the \(i\)th region, expanded at the operating point \(\hat{\kappa}_i\). For any given point \(\kappa = x_0\), the deviation is calculated with respect to the operating point of each region to evaluate each polynomial. The red line in Figure~\ref{fig:regions_poly} represents the deviation of \(x_0\) with respect to the operating point of an arbitrary region, \(\delta \kappa_m\).

{For any arbitrary parameter \(\kappa = x_0\),} the initial states of a family member \(\mathbf{x}\) are evaluated as follows:
\begin{equation}
    \mathbf{x}\left(\kappa = x_0\right) = \sum_{i = 1}^M a_i\left(\kappa\right) \mathcal{P}_i^{\hat{\kappa}_i}\left(\delta \kappa_i\right)
\end{equation}
where \(\delta\kappa_i = \kappa - \hat{\kappa}_i\), as visualized in Figure~\ref{fig:regions_poly}, and \(a_i\left(\kappa\right)\) is the activation function of \(\mathcal{P}_i\), defined as follows:
\begin{equation}
    a_i\left(\kappa\right) = \begin{cases}
    1 & \text{if } \kappa_{il} \leq \kappa < \kappa_{iu}, \\
    0 & \text{otherwise}
\end{cases}
\end{equation}

\subsection{Local Polynomial Regression Model}
The local PRM uses the parameter \(\kappa\) of the designed member as the operating point for the polynomial. This approach yields a single polynomial, which is more computationally efficient than the global method and provides a more precise computation of members that are close to the designed member (i.e., when \(\delta \kappa \ll 1\)). However, this model becomes less effective for larger values of \(\delta \kappa\). The initial states of the member are obtained using this model as follows:
\begin{equation}
    \mathbf{x}\left(\kappa\right) =  \mathcal{P}^{\hat{\kappa}_d}\left(\delta \kappa\right)
\end{equation}
In some cislunar applications, such as station-keeping, only the neighbors of the designed orbit are required, making this local PRM an optimal choice for these cases.

\section{Periodic Orbit Representation Using Differential Algebra}
{As mentioned earlier, a Polynomial Regression Model (PRM) is developed to represent the initial states of each periodic orbit family as a polynomial function of a parameter \(\kappa\). Using this model, the initial state vector is expressed as a function of \(\delta \kappa\), as follows: \(\mathbf{x}_0\left(\delta\kappa\right) = {}_{N}\mathcal{F}^{\hat{\kappa}}\left(\delta\kappa\right)\). Here, \(\mathcal{F}\) denotes the series of the initial states. Instead of propagating the initial states of the periodic orbits individually, it is more efficient to propagate the series \(\mathcal{F}\). In this approach, the resultant series represents the propagated state vector as a function of \(\delta \kappa\) and can be used to obtain the state vector at different values of \(\kappa\) in a less computationally expensive manner.} The DA technique is employed to represent the propagated state vector \(\mathbf{x}_{ij}\) at a certain time \(t_j\), starting from an initial state vector \(\mathbf{x}_{i0}\), as a series expanded at a given parameter vector \(\hat{\kappa}\) to \(N\) order, where each parameter vector \(\mathbf{\kappa}_i\) can be mapped to a certain initial state vector \(\mathbf{x}_{i0}\), as follows:
\begin{equation}
    \begin{aligned}
        & \mathbf{x}_{i0}\left(\delta\mathbf{\kappa}_i\right) = {}_{N}\mathcal{M}^{\hat{\kappa}}_{t_0 \rightarrow t_0}\left(\delta\mathbf{\kappa}_i\right) \\
        & \mathbf{x}_{ij}\left(\delta\mathbf{\kappa}_i\right)
        = {}_{N}\mathcal{M}^{\hat{\kappa}}_{t_0 \rightarrow t_j}\left(\delta\mathbf{\kappa}_i\right)
    \end{aligned}
    \label{eq:mapping_time}
\end{equation}
where \(\mathcal{M}\) denotes the STPM. This STPM is obtained by propagating the CR3BP as a function of the PRM of the periodic orbit family using the fourth-order Runge-Kutta scheme. Therefore, the STPM returns a {mapped} state vector that must be in the subspace of the periodic orbit family for any arbitrary \(\kappa\) is the domain. {This mapping is initially performed at discretized times; however, it is used to obtain the states at any randomly selected time by interpolating the states of the surrounding points.}

\subsection{{Propagation with Respect to Normalized Time}}
{Equation~\eqref{eq:mapping_time} represents the mapping of the propagated states at a given time. In this context, any deviation in \(\kappa\) might lead to a significant deviation in \(\mathbf{x}\) due to the different time periods of the deviated members. This large deviation increases the control effort required to transfer to the deviated member. To minimize the deviation of the propagated states, the STPM can be computed at a certain dimensionless normalized time \(\eta\), where \(\eta = t/{T_{p}}\). In this case, the period \(T_p\) is a function of \(\delta \kappa\) as follows:} 
\begin{equation}
    T_{pi}\left(\delta\mathbf{\kappa}_i\right) = {}_{N}\mathcal{T}^{\delta\hat{\kappa}}\left(\delta\mathbf{\kappa}_i\right)
\end{equation}
{where \({}_{N}\mathcal{T}^{\delta\hat{\kappa}}\left(\delta\mathbf{\kappa}_i\right)\) is the map function of the time period.} 
In this approach, the number of time steps is fixed, while the sampling time is variable and is determined as a ratio of the time period function, as follows:
\begin{equation}
    T_s\left(\delta\kappa\right) = \dfrac{1}{N_s} {}_N\mathcal{T}^{\hat{\kappa}}\left(\delta \kappa\right)
\end{equation}
where \(N_s\) denotes the fixed number of time steps per period. In this case, the resultant state vector at any normalized time can be obtained by adjusting the number of time steps (e.g., \(N_s|_{10\%} = \frac{1}{10}N_s\)). {In this case, the mapping from time to normalized time is performed as follow:}
\begin{equation}
    \eta\left(\delta\mathbf{\kappa}_i\right) = \dfrac{t}{{}_N\mathcal{T}^{\hat{\kappa}}\left(\delta \kappa\right)}
\end{equation}
The propagated states at a certain \(\eta\) is represented as follows:
\begin{equation}
    \begin{aligned}
        & \mathbf{x}_{i0}\left(\delta\mathbf{\kappa}_i\right) = {}_{N}\overline{\mathcal{M}}^{\hat{\kappa}}_{\eta_{0} \rightarrow \eta_{0}}\left(\delta\mathbf{\kappa}_i\right) \\
        & \mathbf{x}_{ij}\left(\delta\mathbf{\kappa}_i\right) =  {}_{N}\mathcal{M}^{\hat{\kappa}}_{t_0 \rightarrow t_j}\left(\delta\mathbf{\kappa}_i, T_s\left(\delta\mathbf{\kappa}_i\right)\right) = {}_{N}\overline{\mathcal{M}}^{\hat{\kappa}}_{\eta_{0} \rightarrow \eta_{j}}\left(\delta\mathbf{\kappa}_i\right)
    \end{aligned}
    \label{eq:mapping_normalizedTime}
\end{equation}

where \(\overline{\mathcal{M}}\) denotes the STPM using normalized time. {Figure~\ref{eq:mapping_normalizedTime} depicts the difference between mapping to time and normalized time, with the horizontal solid lines representing points at the same normalized time \(\eta\), while the dashed lines representing points at the same time \(t\).} 

\begin{figure}[H]
    \centering
    \includegraphics[width=0.4\linewidth]{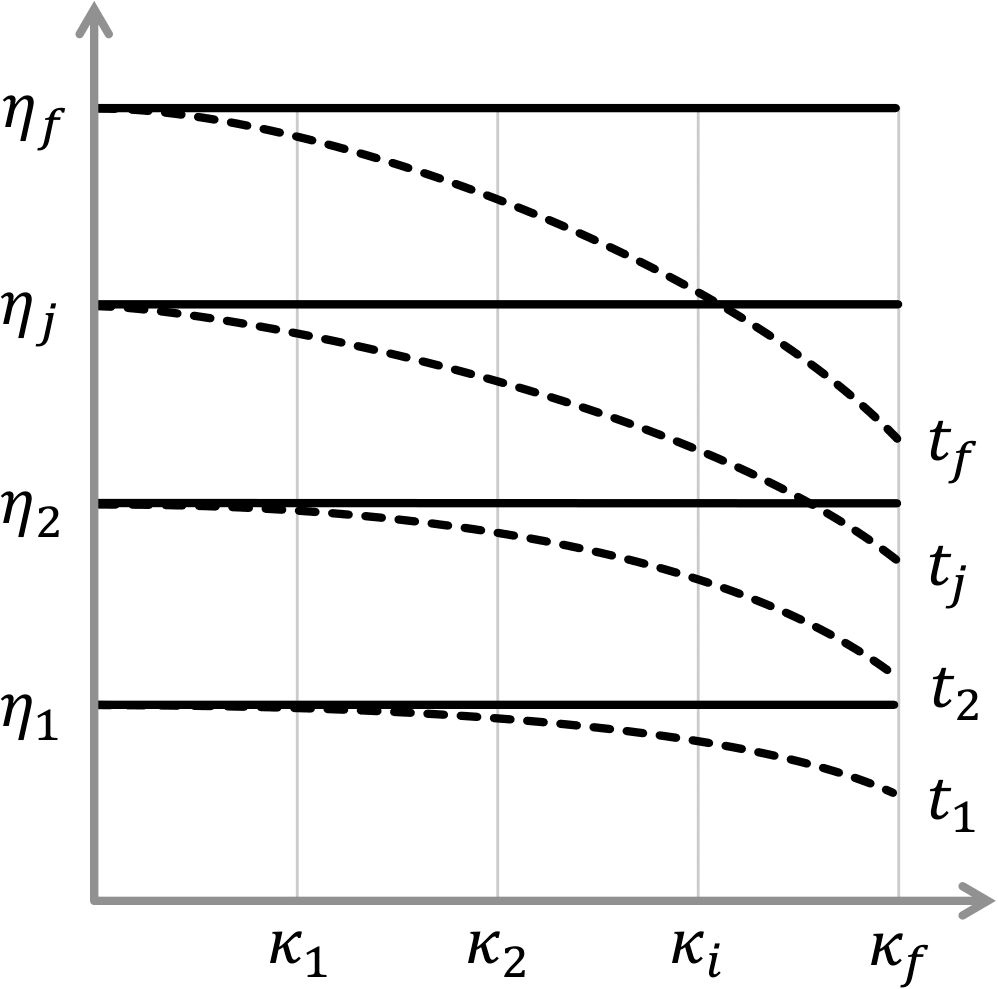}
    \caption{A comparison between mapping to time versus mapping to normalized time.}
    \label{fig:timeVsNormalizedTime}
\end{figure}

\section{Numerical Simulations}
The proposed PRMs and the DA representation of periodic orbit families are verified through a series of numerical simulations.

\subsection{Precision of the Derived  Polynomial Regression Models}
\begin{figure}[H]
    \centering
    \subfigure[Propagation for three orbits.]{%
        \includegraphics[width=0.45\textwidth]{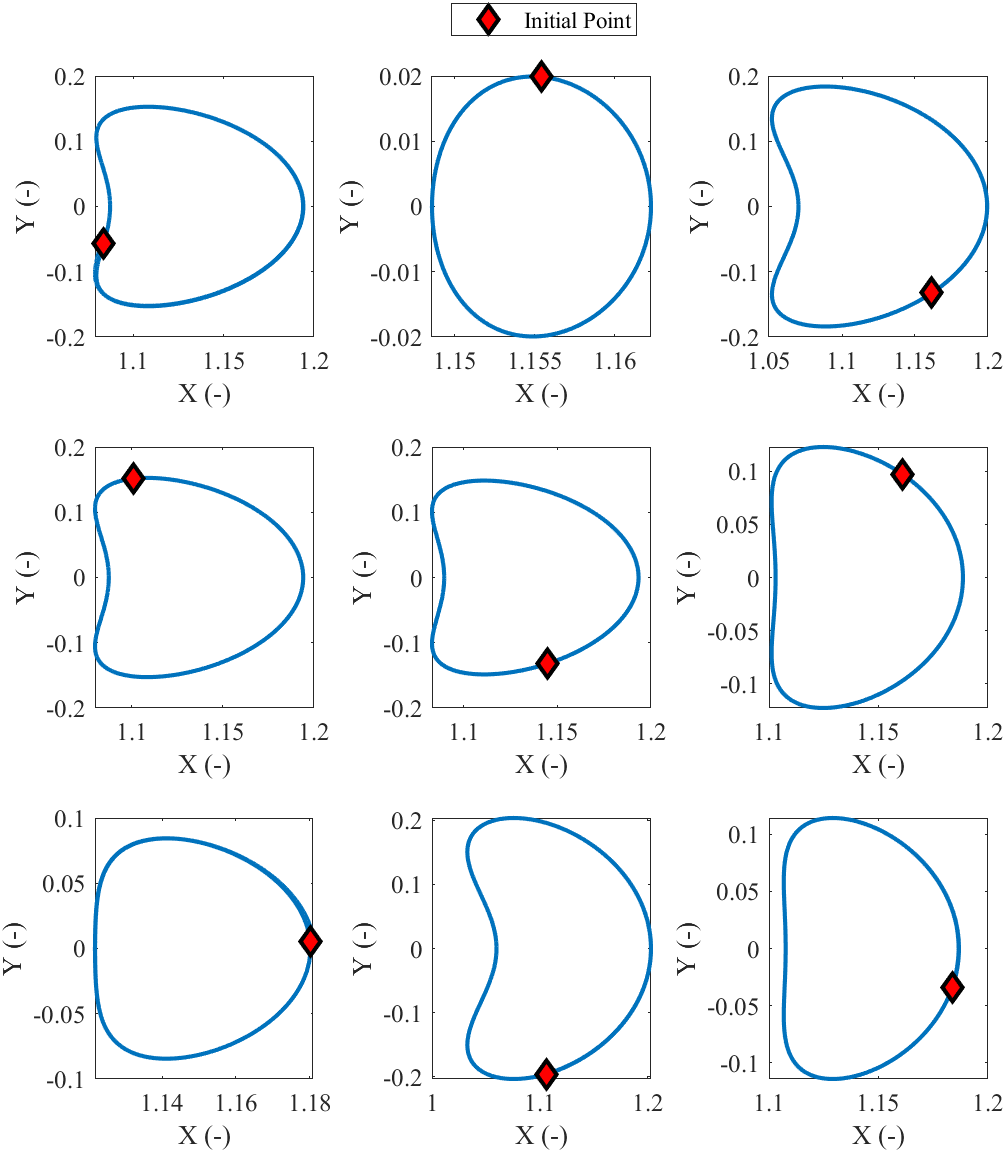}
        \label{fig:propagate_3O}
    }
    \hfill
    \subfigure[Propagation for four orbits.]{%
        \includegraphics[width=0.45\textwidth]{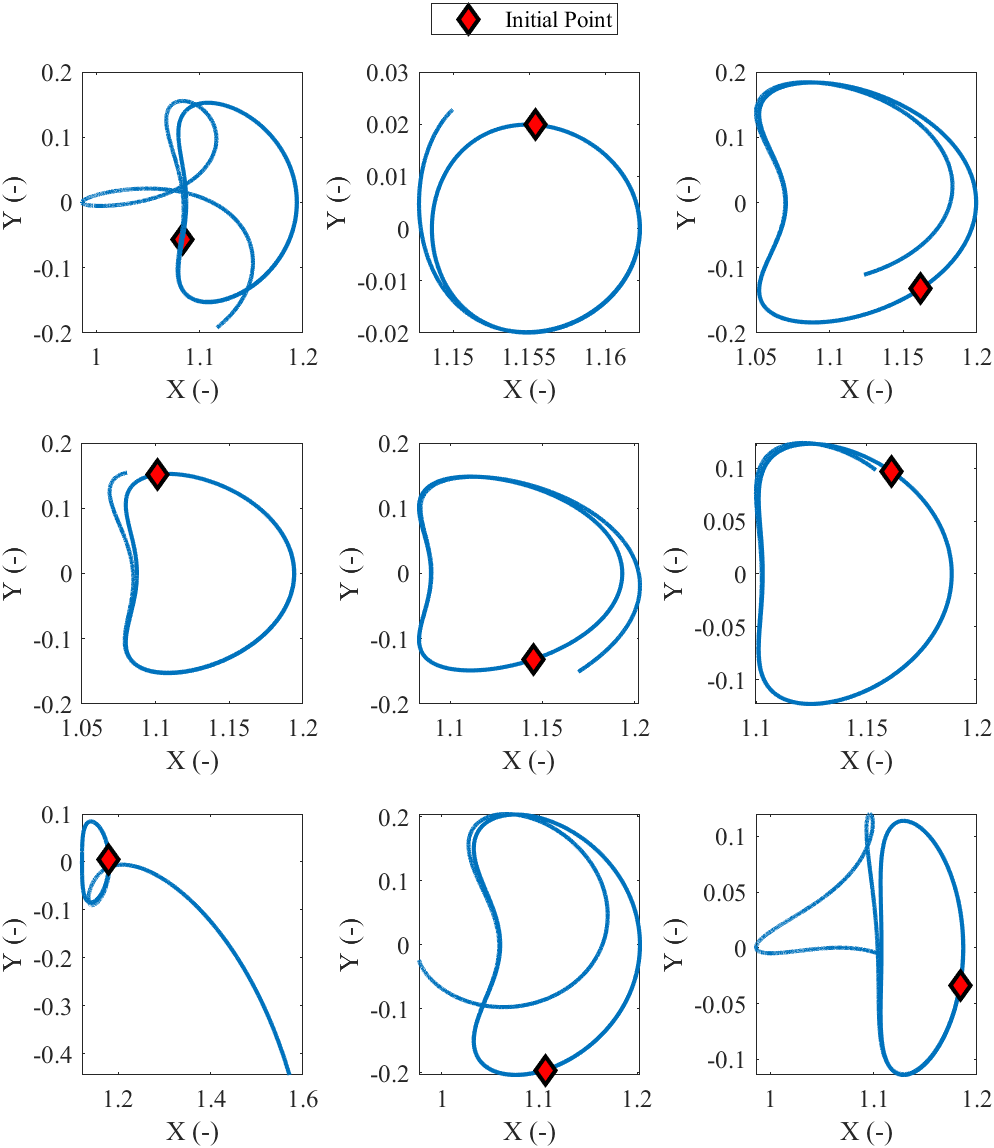}
        \label{fig:propagate_4O}
    }
    \caption{The propagated states using the PRM.}
    \label{fig:propagate}
\end{figure}

\begin{figure}[H]
    \centering
    \includegraphics[width=0.54\linewidth]{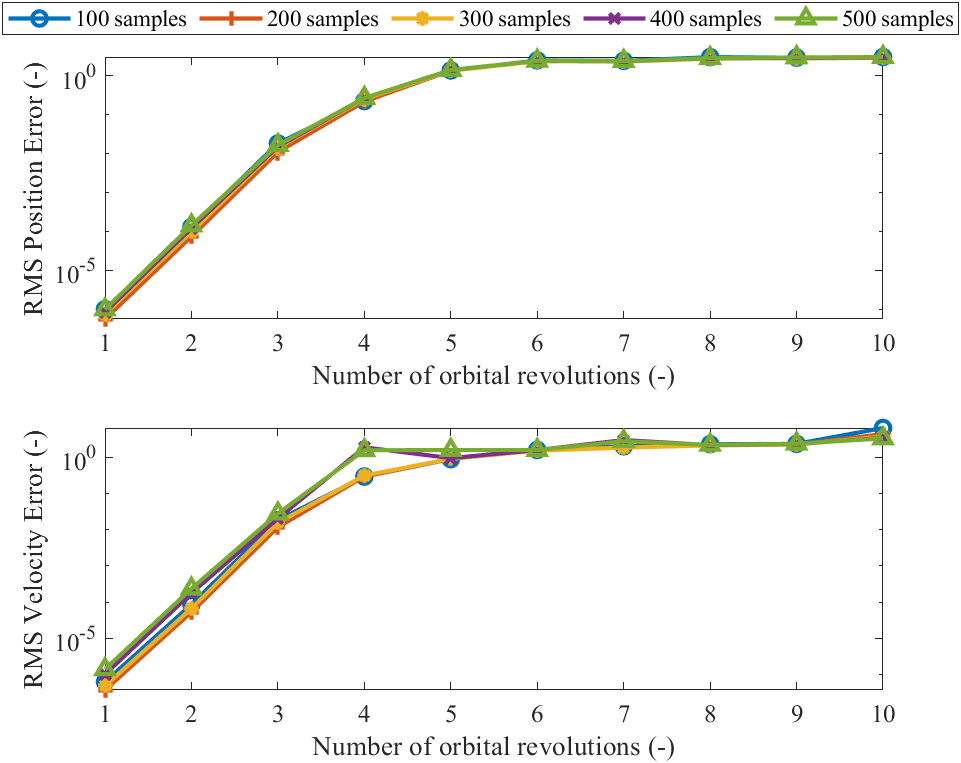}
    \caption{The RMS error of the states that are randomly initialized by the PRM and propagated for multiple orbital revolutions.}
    \label{fig:propagate_err}
\end{figure}

{The purpose of the PRM is to obtain the initial states of any member in the periodic orbit family in a computationally efficient manner. In this numerical simulation, the accuracy of the developed model is assessed by measuring the state error after propagation for multiple orbits.}
To evaluate the accuracy of the proposed global PRM, random states are generated at random times and propagated using the RK4 method for a finite number of orbits. {Both random points and times are generated using a uniform distribution in a MATLAB environment.} For the LOs near \(\mathcal{L}_2\), the domain is divided into {eight regions, with polynomials of order thirty.} Figure~\ref{fig:propagate} shows the propagation of these random initial points over time. Figure~\ref{fig:propagate_3O} demonstrates that the states maintained the periodic orbit over three orbits with only insignificant deviations. However, after three orbits, Figure~\ref{fig:propagate_4O} reveals that most of these states fail to maintain the periodic orbits. The accuracy would vary if the number of regions or the polynomial order were changed. { Figure~\ref{fig:propagate_err} shows the Root Mean Square Error (RMSE) of both position and velocity after different numbers of orbital revolutions. In this analysis, initial states are randomly generated using the developed PRM. The initial position vector is randomly selected, and the corresponding velocity vector is then computed using the PRM to satisfy the periodic orbit conditions. The figure indicates that both position and velocity errors follow the same trend and are of similar magnitude for all number of samples. The errors start from significantly small values and settle to an order of magnitude of one. This analysis demonstrates that the developed PRM can compute accurate initial states for periodic orbits, maintaining a significantly small error even after up to three orbits.}


\subsection{Periodic Orbit Families Using the Derived Polynomial Regression Models}
The developed global PRM is applied to compute the members of the Lyapunov and HO families near the libration points \(\mathcal{L}_1\) and \(\mathcal{L}_2\). {The polynomial approximations of these orbits are obtained, representing the initial states of the periodic orbit as a function of \(x_0\).} Then, these polynomials are propagated over time using the RK4 method within a DA framework. To validate the method's accuracy and demonstrate its effectiveness, the polynomials are propagated over various time and normalized time intervals. This approach provides a comprehensive view of the orbits' trajectories and their stability characteristics. The results, showcasing the generated polynomials, are presented in this section, highlighting the efficiency and precision of the global and local PRMs in orbit computation and propagation near the libration points.

\subsection{Lyapunov Orbits Family}
\begin{figure}[H]
    \centering
   \includegraphics[width=0.9\linewidth]{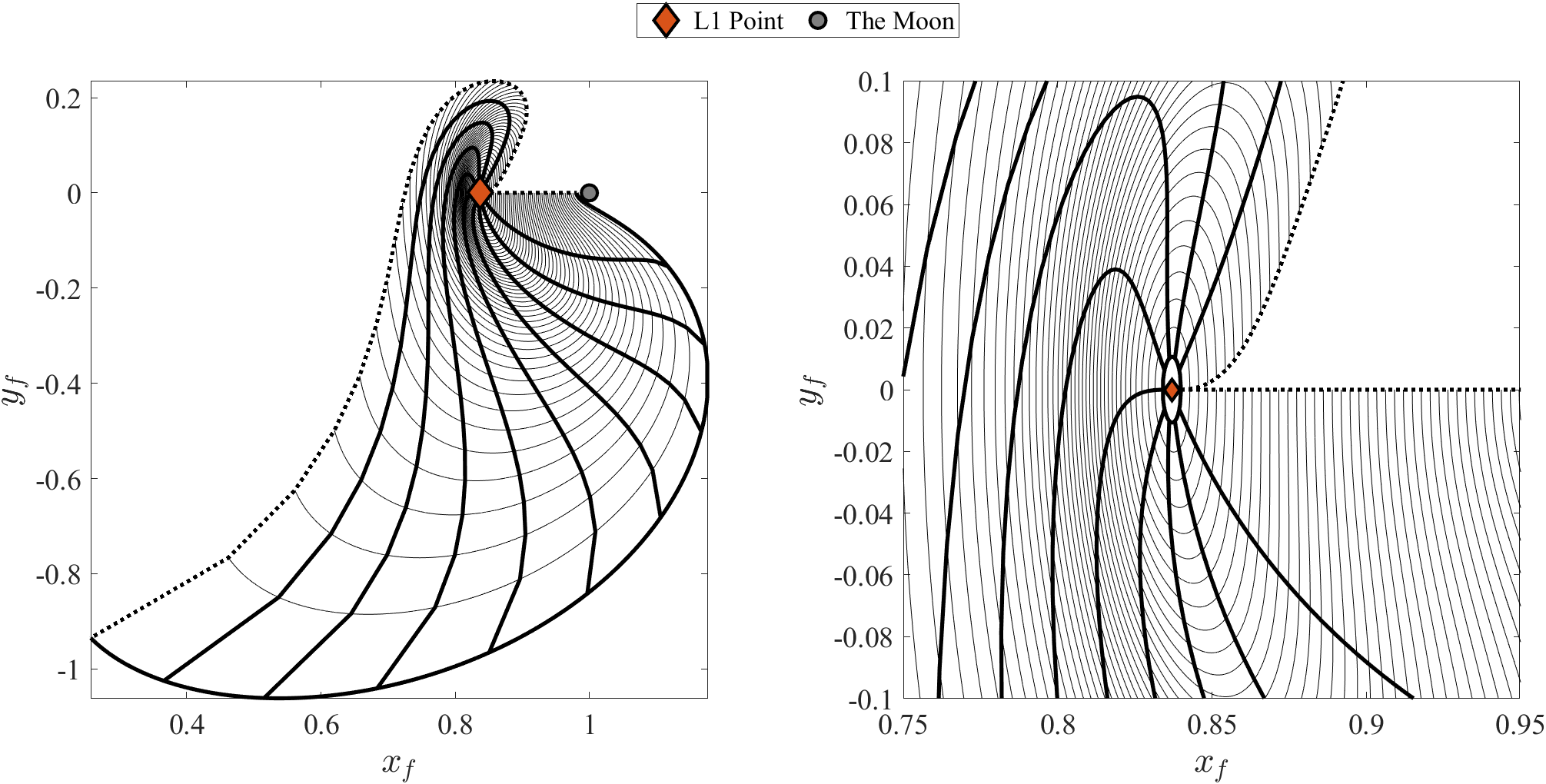}
    \caption{The propagated states to different times using the PRM.}
    \label{fig:allOrbits_t_L1_Lyap}
\end{figure}

\begin{figure}[H]
    \centering
    \includegraphics[width=0.45\linewidth]{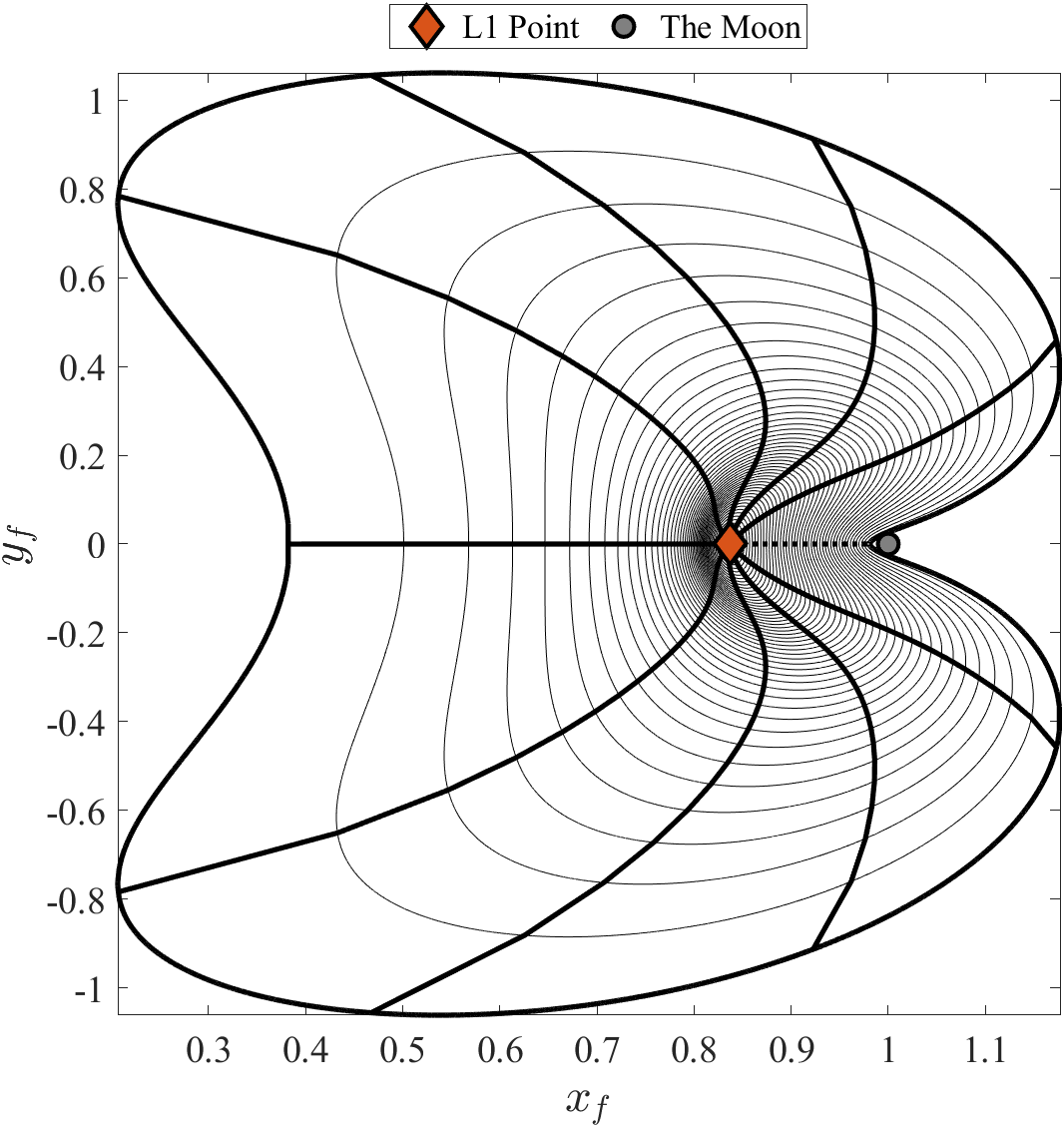}
    \caption{The propagated states to different normalized times using the PRM.}
    \label{fig:allOrbits_tau_L1_Lyap}
\end{figure}

\begin{figure}[H]
    \centering
    \includegraphics[width=\linewidth]{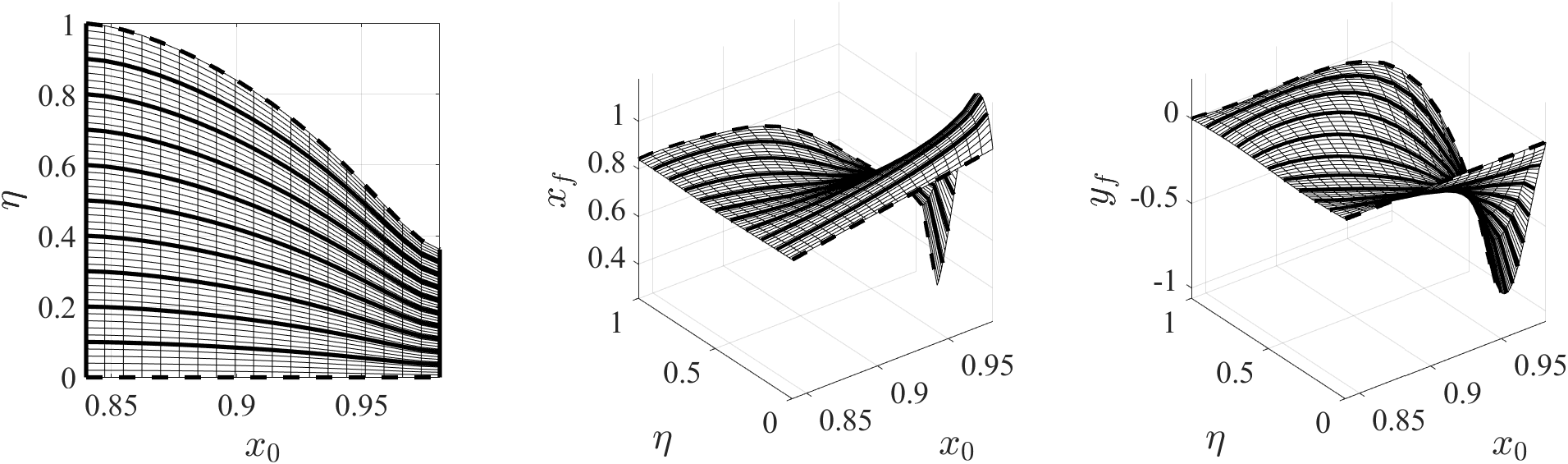}
    \caption{The variation of the position vector components with \(x_0\) and the normalized time.}
    \label{fig:Az_vs_tau_state_L1_Lyap}
\end{figure}

\begin{figure}[H]
    \centering
    \includegraphics[width=\linewidth]{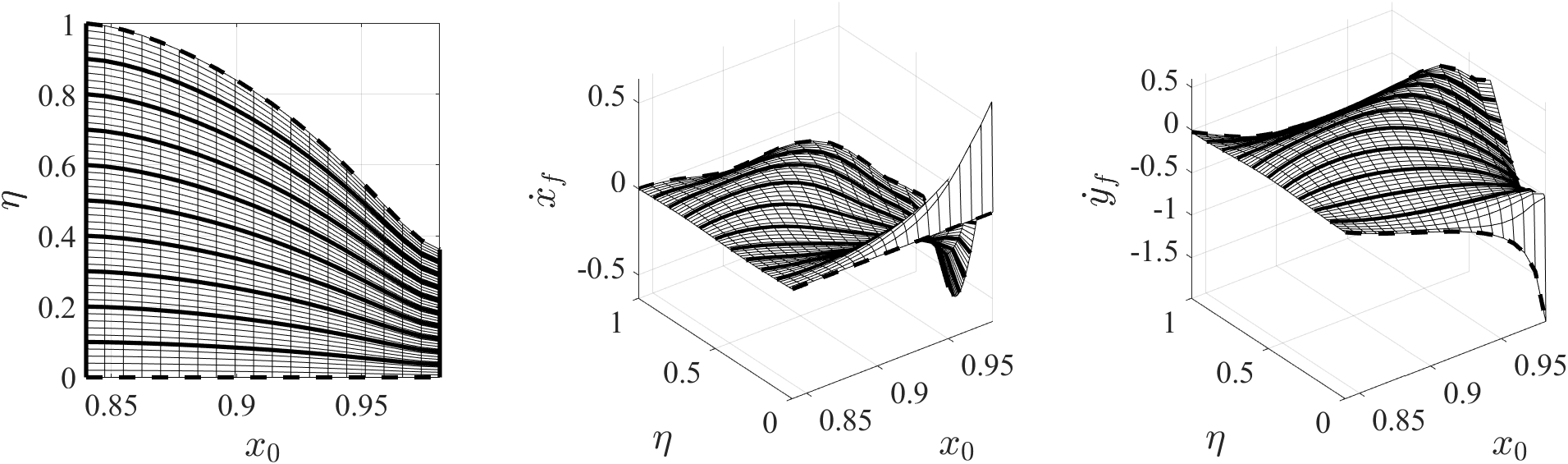}
    \caption{The variation of the velocity vector components with \(x_0\) and the normalized time.}
    \label{fig:Az_vs_tau_state_dot_L1_Lyap}
\end{figure}

In LOs, the initial condition {on the Earth-Moon rotating frame} is governed by two states: \(x_0\) and \(\dot{y}_0\), as they are planar orbits. Since the predefined parameter of the PRM is \(x_0\), the state \(\dot{y}_0\) is expressed as a polynomial. However, after propagating the states to a certain time, each of the four planar states will have a different polynomial expressed as a function of \(\delta\kappa\) or \(\delta x_0\). Figure~\ref{fig:allOrbits_t_L1_Lyap} shows the states of the LOs near \(\mathcal{L}_1\) after propagating the initial states to ten different times. Each solid line represents the locus of points that share the same time. Figure~\ref{fig:allOrbits_tau_L1_Lyap} displays the propagation of the states to ten different normalized times, with each solid line representing the locus of points that share the same normalized time. {It is worth noting that normalized time propagation covers the entire domain of the family. Additionally, points with the same normalized time are closer to each other than points with the same time, especially in long-term propagation. This demonstrates the superiority of normalized time propagation over time propagation in proximity operations.} Figure~\ref{fig:Az_vs_tau_state_L1_Lyap} illustrates the variation of the position states with \(x_0\) and \(\eta\), and also demonstrates how \(\eta\) varies with \(x_0\), with each solid line representing normalized times of the family at a given time. Figure~\ref{fig:Az_vs_tau_state_dot_L1_Lyap} shows the variation of the velocity states with \(x_0\) and \(\eta\) in a similar manner. Here, \(x_0\) den
\subsection{Halo Orbit Families}
\begin{figure}[H]
    \centering
    \subfigure[Propagation with time.]{%
        \includegraphics[width=0.45\textwidth]{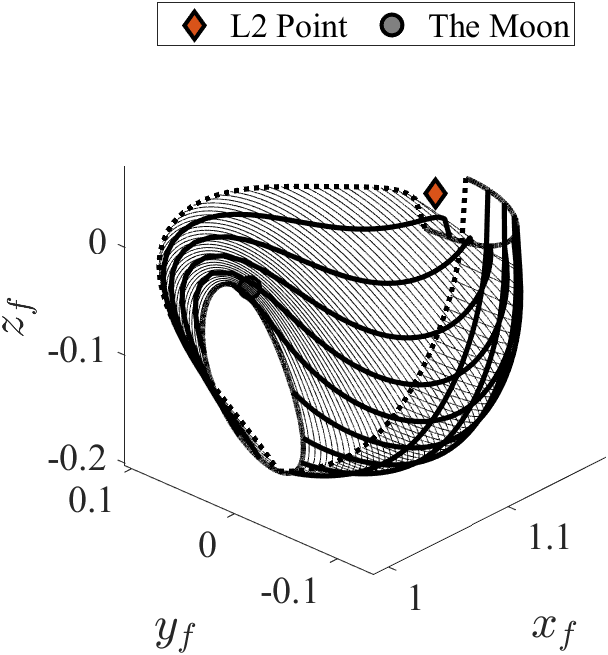}
        \label{fig:allOrbits_t_L2_Halo}
    }
    \hfill
    \subfigure[Propagation with normalized time.]{%
        \includegraphics[width=0.45\textwidth]{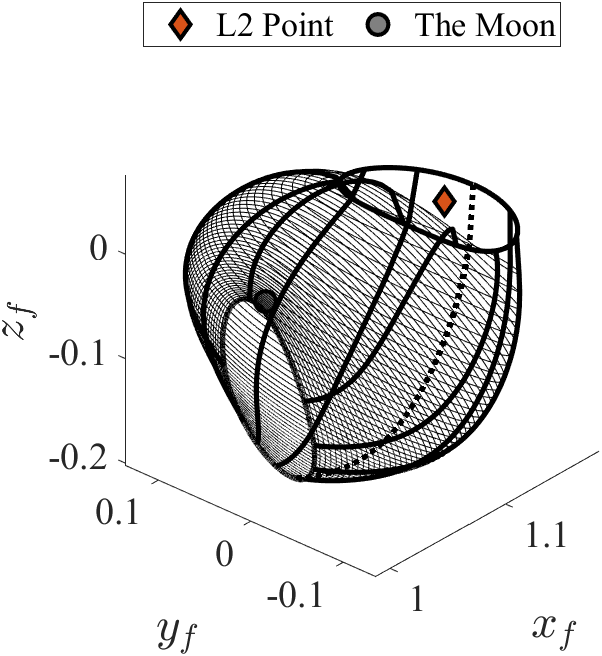}
    }
    \caption{The propagated states to different times and normalized times using the PRM.}
    \label{fig:allOrbits_tau_L2_Halo}
\end{figure}

\begin{figure}[H]
    \centering
    \includegraphics[width=\linewidth]{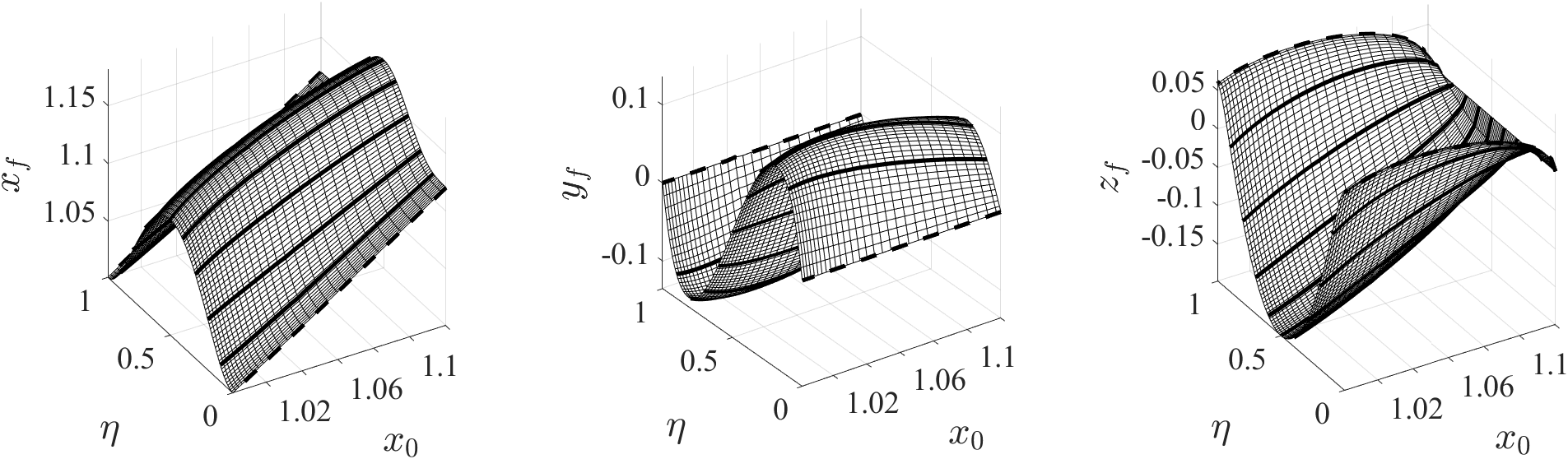}
    \caption{The variation of the position vector components with \(x_0\) and the normalized time.}
    \label{fig:Az_vs_tau_state_L2_Halo}
\end{figure}

\begin{figure}[H]
    \centering
    \includegraphics[width=\linewidth]{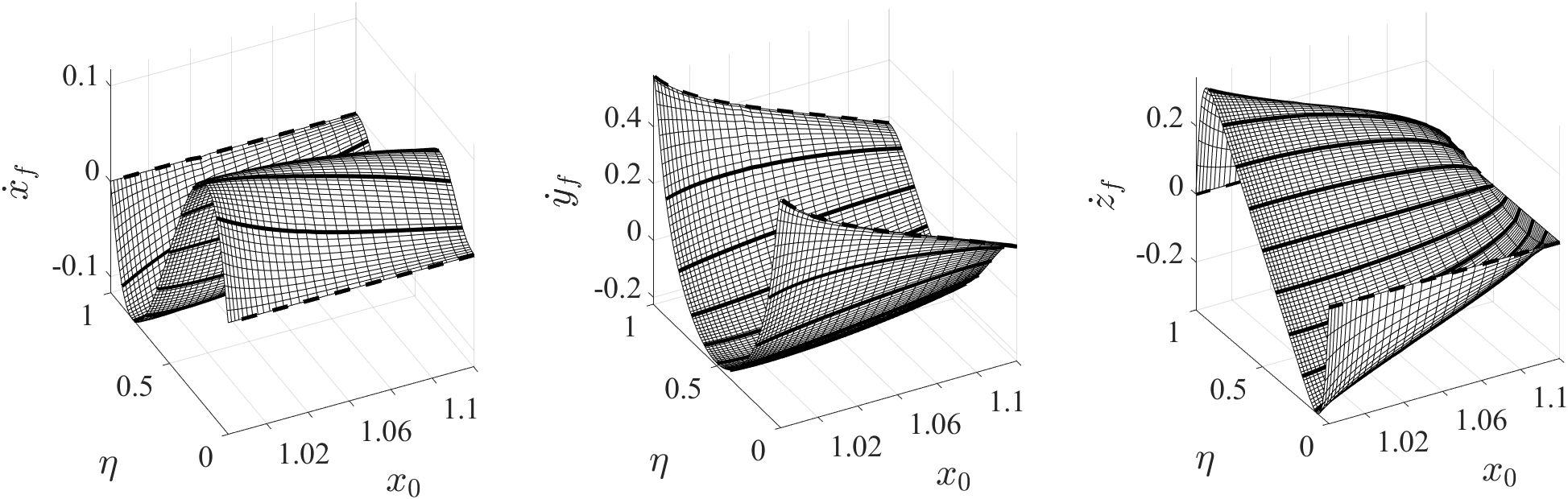}
    \caption{The variation of the velocity vector components with \(x_0\) and the normalized time.}
    \label{fig:Az_vs_tau_state_dot_L2_Halo}
\end{figure}

In HOs, the initial conditions are governed by three states: \(x_0\), \(z_0\), and \(\dot{y}_0\), as these are 3-D orbits. The states \(z_0\) and \(\dot{y}_0\) are defined as functions of \(\delta x_0\). After propagation, each of the six states will have a polynomial representation at each time or normalized time. Figure~\ref{fig:allOrbits_t_L2_Halo} shows the states of the HOs near \(\mathcal{L}_2\) after propagating the initial states to ten different times, with each solid line representing the locus of points that share the same time. Figure~\ref{fig:allOrbits_tau_L2_Halo} illustrates the propagation of the states to ten different normalized times, with each solid line representing the locus of points that share the same normalized time. Figure~\ref{fig:Az_vs_tau_state_L2_Halo} shows the variation of the position states with \(x_0\) and \(\eta\), while Figure~\ref{fig:Az_vs_tau_state_dot_L2_Halo} shows the variation of the velocity states.


\subsection{Comparative Analysis of Global and Local Polynomial Regression Models}
\begin{figure}[H]
    \centering
    \includegraphics[width=0.6\linewidth]{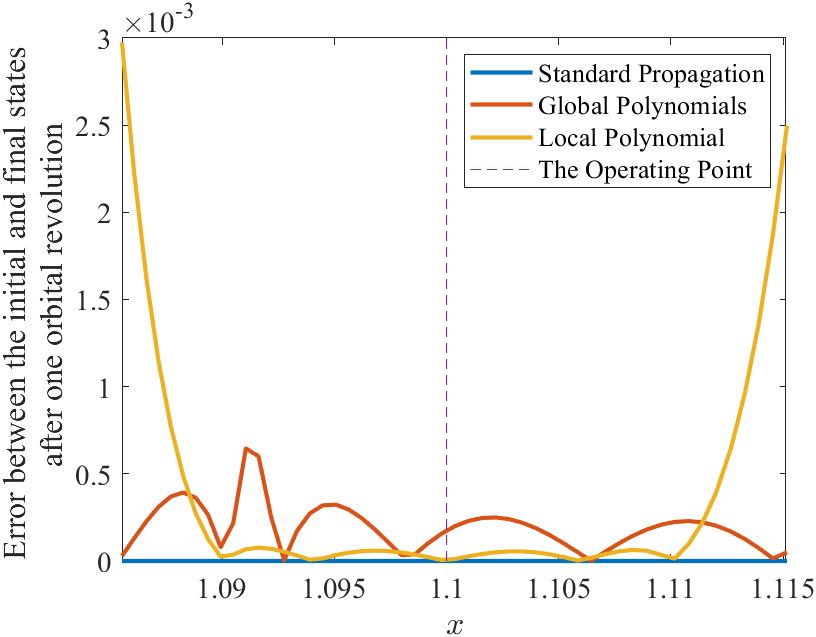}
    \caption{A comparison between the precision of propagating the states using global and local PRMs.}
    \label{fig:globalVsLocal}
\end{figure}

The strength of the global PRM lies in its ability to cover the entire domain of the orbit family. However, locating the orbit in the correct region can be computationally expensive for some applications. Conversely, if an application only requires information about the neighborhood of a specific design orbit rather than the entire family, the local PRM becomes the optimal choice, {since it uses a single polynomial, making it computationally efficient.} Figure~\ref{fig:globalVsLocal} compares the accuracy of the global and local PRMs, represented by the error between the final and initial states after one time period. An arbitrary point is selected for comparison, and its neighboring points are propagated using the three methods. The comparison shows that the global model maintains consistent accuracy across the domain, whereas the local models exhibit better accuracy near the designed operating point. However, as expected, the local model fails to provide accurate results when the states deviate significantly from the designed operating point. It is worth noting that the global PRM requires evaluating multiple polynomials to maintain consistent error, while the local PRM requires evaluating only a single polynomial.

\subsection{Proportional-Derivative Controller}
A PD control law is designed to verify the applicability of the proposed models in real-world missions. The controller gains are manually tuned until the system is stabilized. The objective of the controller is to transfer the satellite to one of the nearest periodic orbits while keeping it within the same family. The control law is implemented as follows:
\begin{enumerate}
    \item The position of the satellite is measured to identify the closest periodic orbit, as determined by the proposed model.
    \item The nearest periodic orbit parameter, \(\kappa_0\), and its normalized time, \(\eta_0\), are determined using the normalized time mapping, as depicted in Figure~\ref{fig:allOrbits_tau_L1_Lyap} and Figure~\ref{fig:allOrbits_tau_L2_Halo}.
    \item The normalized transfer time, \(\eta_{t}\), is selected; then, the reference state is obtained as follows: $$\mathbf{x}_{r} = {}_{N}\overline{\mathcal{M}}^{\hat{\kappa}}_{\eta_{0} \rightarrow \eta_{0} + \eta_{t}}\left(\delta\mathbf{\kappa}_0\right)$$
    \item The thrust is evaluated as follows: $\mathbf{u} = - \left[\mathbf{K}\right] \times \left(\mathbf{x} - \mathbf{x}_{r}\right)$,
    where $\left[\mathbf{K}\right] = \text{Diag}\left(k_p\;\mathbf{I_3},\;\;k_d\;\mathbf{I_3}\right)$. Here, \(\mathbf{I_3}\) denotes the identity matrix, and \(k_p\) and \(k_d\) are the controller gains.
    \item The equivalent transfer time is computed as follows: $t_t = \eta_t\times{}_N\mathcal{T}^{\hat{\kappa}}\left(\delta \kappa_0\right)$.
    \item The velocity impulse \(\Delta \mathbf{V}\) is evaluated as follows: $\Delta \mathbf{V} = \mathbf{u}\; t_t$, which is a valid approximation for small \(t_t\).
    \item These steps are repeated frequently after each transfer.
\end{enumerate}

\begin{figure}[H]
    \centering
    \subfigure[The trajectory.]{%
        \includegraphics[width=0.44\textwidth]{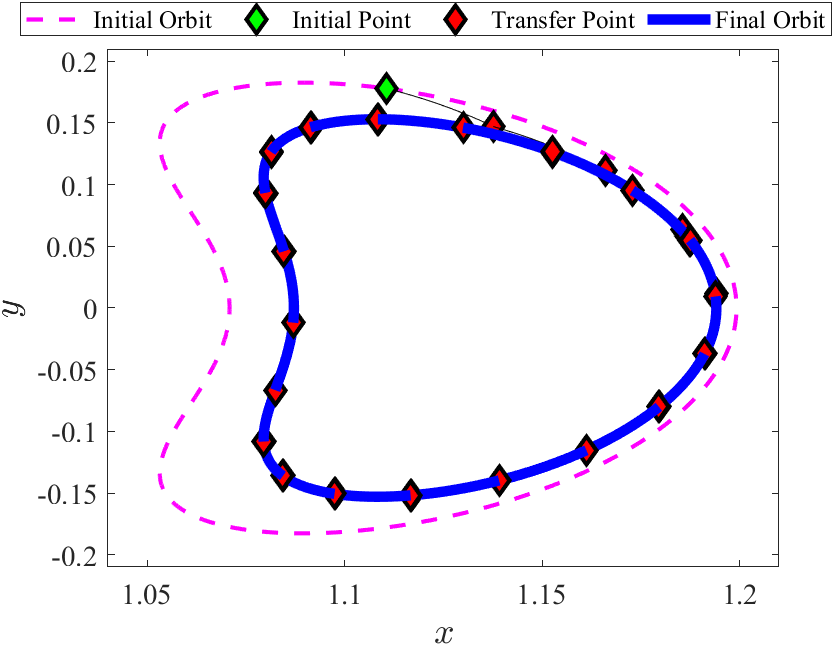}
        \label{fig:controlOrbits_traj}
    }
    \hfill
    \subfigure[The control efforts.]{%
        \includegraphics[width=0.44\textwidth]{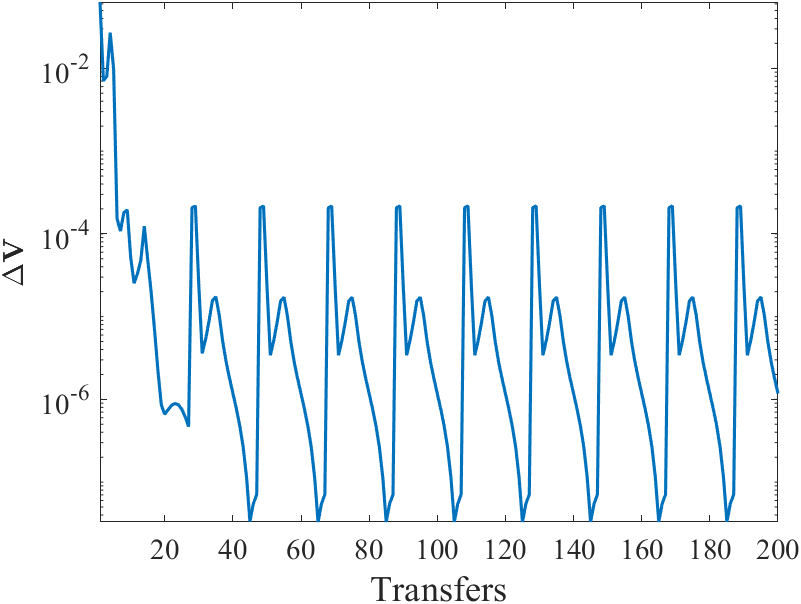}.
        \label{fig:controlOrbits_delV}
    }
    \caption{The trajectory and the time history of control efforts obtained using the PD control law.}
    \label{fig:controlOrbit}
\end{figure}

The control law is tested on the LO family near \(\mathcal{L}_2\), where a random periodic orbit is selected, and the satellite starts from a random point on that orbit with a disturbed velocity. The control law is then implemented to retain the satellite within the family by performing multiple transfers. The transfer time is set to \(\eta_{t} = 0.05\). Figure~\ref{fig:controlOrbits_traj} shows the transfers over ten orbits until the satellite is retained in a specific periodic orbit. Figure~\ref{fig:controlOrbits_delV} illustrates the time history of the velocity impulses, indicating that the satellite rapidly converges to a steady-state periodic orbit within two revolutions. Although the final periodic orbit differs from the initial one, it maintains the same ground track on the Moon. These results emphasize the potential of the proposed method for various cislunar applications, including low-energy transfers.

\begin{figure}[H]
    \centering
    \includegraphics[width=0.8\linewidth]{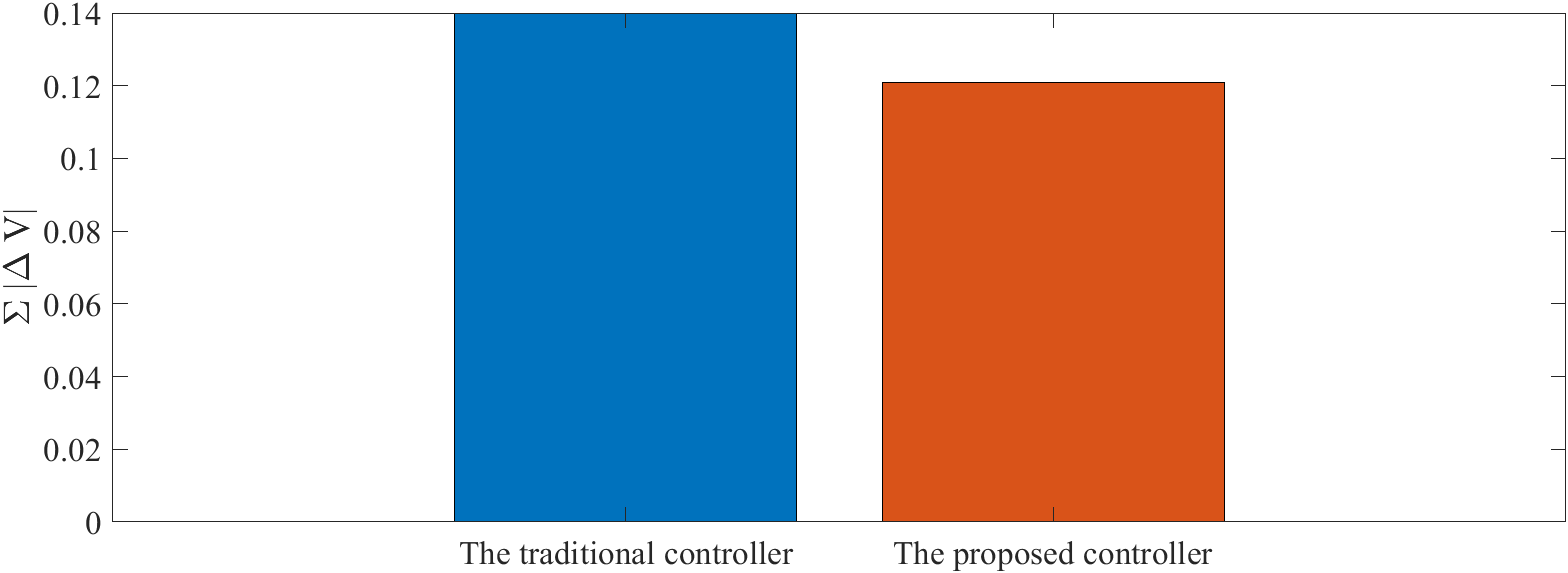}
    \caption{A comparison of control effort between the proposed controller and the traditional controller.}
    \label{fig:TradVsNewRand}
\end{figure}

The proposed control law is compared to a traditional controller, which tracks the states of a predetermined orbit without using the STPM. As previously mentioned, the proposed method converges to a random orbit within the family. To efficiently compare the two methods, the same simulation parameters are used, and the final orbit obtained from the proposed method is set as the target orbit for the traditional method. Figure~\ref{fig:TradVsNewRand} presents the control efforts of both methods, demonstrating a significant reduction of \(13.54\%\) in the impulses required by the proposed method.

\section{Conclusions}
In conclusion, this study explores the dynamics of motion in periodic orbits near libration points in cislunar space using the DA framework. The CR3BP model is employed to describe the motion in this environment. Initial states of the periodic orbits are numerically generated using a differential correction scheme that leverages an analytical solution, while the members of each orbit family are computed using the PAC method. These computed members are then used to fit PRMs for the orbit families, with the initial states expressed as functions of predefined parameters. These regression models are incorporated into the DA framework to evaluate propagated states at a given time as functions of deviations in these predefined parameters. The accuracy of the computed states at various times is assessed, and the execution time for computing these states is compared with traditional propagation methods using the Runge-Kutta method. The analysis demonstrates the effectiveness of using DA for representing motion in periodic orbits in cislunar space. Additionally, it shows significantly reduced control effort compared to the
traditional tracking control law.

\section{Acknowledgement}
The authors wish to thank Erin Ashley and Batu Candan for their help in reviewing and improving this paper.

\bibliographystyle{ieeetr}
\bibliography{references}

\end{document}